

\documentstyle{amsppt}
\magnification=\magstep 1
\TagsOnRight
\NoBlackBoxes
\leftheadtext{V\.G\. Pestov}
\rightheadtext{Odd coordinates}
\def\norm #1{{\left\Vert\,#1\,\right\Vert}}

\def\R {{\Bbb R}}
\def\C {{\Bbb C}}

\def\N{{\Bbb N}}
\def\e{{\epsilon}}

\def\G {{\bigwedge l_1}}
\def\H{{\widehat{\bigwedge l_1}}}
\def\X{{\frak X}}
\def\s{{\Sigma (\H)}}

\topmatter
\title
An analytic structure emerging in presence of
infinitely many odd coordinates
\endtitle
\author
Vladimir G\. Pestov
\endauthor
\affil
Department of Mathematics, Victoria University of Wellington,
P.O. Box 600, Wellington, New Zealand \\ \\
vladimir.pestov$\@$vuw.ac.nz
\endaffil
\abstract{We show that the spectrum of the locally convex
nonstandard
hull of an infinite dimensional Grassmann algebra
contains a nontrivial analytic part.}

\endabstract
\endsubjclass
\endkeywords
\endtopmatter
\document

\smallpagebreak
\heading
\S 1. Introduction
\endheading
\smallpagebreak

The present note contributes to the program of
featuring even geometry as a ``collective effect in
infinite-dimensional odd geometry,'' as
suggested by Manin \cite{14}.
In our earlier work \cite{19, 22} we proposed to apply the construction of
nonstandard hull \cite{12} to a normed Grassmann algebra
$\bigwedge (\nu)$ with an infinitely large number of odd generators,
$\nu$, in order to construct the superspace $\X$ limit of
finite dimensional purely
odd superspaces as
a locally ringed superspace associated to the nonstandard hull algebra
$\widehat{\bigwedge (\nu)}$ in a univeral way.
We have demonstrated that the resulting superspace
$\X$ is infinite-dimensional
both in odd sector and, remarkably, in even sector as well, in
the sense that
the spectrum of the hull algebra ${\widehat{\bigwedge (\nu)}}$,
which serves as the underlying compact space for $\X$,
is infinite-dimensional.

It seems, however, that the  nonstandard hull of
the normed
Grassmann algebra admits practically no derivations
other than nilpotent, which makes it hardly possible to
develop a far-reaching
version of ``even'' analysis and geometry over it.
In this note we propose to consider instead a nonstandard hull
formed with respect to a certain natural
non-normable locally convex topology
on an infinite-dimensional Grassmann algebra.
Namely, we start with the (completed) exterior algebra
$\G$ over a countably infinite
set $\Xi$ of odd generators which has the strongest
locally convex algebra topology making $\Xi$ into a bounded set.
The spectrum, $\s$,
of the resulting locally convex nonstandard hull $\H$
carries a nontrivial analytic structure, namely,
there exists a homeomorphism $\Phi$ from the unit disc $D\subset \C$
to $\s$ such that for every element
$a\in\H$ the composition $\hat a\circ\Phi$ is a holomorphic
function from $D$ to $\C$.
The presence of such analytic parts in the spectrum of
a locally convex algebra is known to be uncommon
(cf. \cite{29}).
The result indicates that infinite
dimensional purely odd superspaces
may in a sense ``generate'' the complex analyticity.

For an elaborate discussion of the present circle of ideas
the reader is referred to our paper \cite{22}.
We follow the ``classical'' Robinson's approach to
nonstandard analysis \cite{26, 30}, because apparently no other known
formalization of nonstandard analysis (such as Internal Set
Theory) provides adequate means
to deal with nonstandard hulls.
Either of the books \cite{2, 15} serves as a
self-contained introduction to supergeometry in the form we need.
Terminology and notation from theory of topological
algebras, locally convex spaces,
infinite-dimensional holomorphy, and Banach-Lie theory follow
\cite{29}, \cite{27}, \cite{5}, and \cite{3}, respectively.
Throughout the paper, $\C$ is the basic field.

\smallpagebreak
\heading
\S 2. Free locally convex graded commutative algebras
\endheading
\smallpagebreak

A rich collection of ``grassmannian'' algebras can be obtained
by invoking the categorial concept of a universal arrow to
forgetful functor \cite{13}.
Here we will present just one of them, and we refer the reader for
further examples of the kind to \cite{9, 20, 21, 23, 24}.
The present particular
construct seems to be new, although numerous similar
algebras, known for quite a while in the folklore, are now
surfacing in mathematical physics every now and then.
\footnote{The author has benefited from
discussing this topic
with J. Baez, U. Bruzzo, Z. Jaskolski, T. Loring, and N.C. Phillips.}

Recall that a $\Bbb Z_2$-graded algebra $\Lambda$
(we will term such
algebras  just {\it graded}) is {\it graded commutative}
if for all
$x,y\in\Lambda$ one has
$xy = (-1)^{\tilde x\tilde y}yx,~~ x,y \in \Lambda^0 \cup
\Lambda^1 $,
where $\tilde x\in \Bbb Z_2$ is the {\it parity of} $x$,
such that $x\in\Lambda^{\tilde x}$.

\proclaim{2.1. Theorem}
Let $E$ be a locally convex space. There exist
a complete locally convex graded commutative unital algebra
$\bigwedge E$ and a continuous linear operator
$i:E\to\bigwedge E^1$
such that
every continuous linear operator $f$ from $E$ to the odd
part $\Lambda^1$ of a complete locally convex
graded commutative unital algebra $\Lambda$
gives rise in a unique way to a continuous unital graded algebra
homomorphism $\hat f: \bigwedge E \to \Lambda$ with
$f=\hat f\circ i$.
\qed\endproclaim

The proof employs certain pretty standard
devices
(cf. \cite{20, 21, 23, 24}).
What really matters, is the following result
which details the structure
of the free locally convex graded commutative algebra
$\bigwedge E$.

\proclaim{2.2. Theorem}
Let $E$ be a locally convex space. The following are true.
\item{(i)}
The mapping $i:E\to\bigwedge E^1$ is an embedding of locally convex
spaces;
\item{(ii)} the algebra generated by $i(E)$ in
$\bigwedge E$ is the exterior algebra over $E$, and it is
dense in $\bigwedge E$;
\item{(iii)} the
completed $n$-th symmetric power of $i(E)$ in  $\bigwedge E$,
$\bigwedge^n(E)$, is the locally convex quotient space
of the $n$-th projective tensor power $E^{\hat\otimes n}$
under the antisymmetrization map;
\item{(iv)} the algebra $\bigwedge E$ is the locally convex direct sum
of subspaces $\bigwedge^n(E)$, $n\in\N$.
\endproclaim

\demo{Proof}
The only nontrivial point is the verification of the continuity of
multiplication in an algebra constructed as
the locally convex direct sum
of subspaces $\bigwedge^n(E)$, $n\in\N$. This is done with the
help of the fact that on $\bigwedge E$ three topologies coincide:
the topology of the locally convex direct sum, the topology of
the direct limit of the
spaces $\oplus_{i=0}^n\bigwedge^i(E),~i\in\N$
in the category of Tychonoff spaces and the
box product topology (see \cite{27}).
\qed\enddemo

We will mention just two examples. If $E=\C^n$ is finite dimensional
then the algebra $\bigwedge E$ is nothing else but the familiar
Grassmann algebra with $n$ generators.
The case where $E=\C^\infty \equiv \varinjlim\C^n$
was actually studied in \cite{10}.
For other related examples, see \cite{24}.

\proclaim{2.3. Theorem} Let $\Gamma$ be a set.
Then there exists a locally convex unital graded commutative
algebra $\Lambda$ topologically
generated by its bounded set $\Gamma$ of odd elements,
such that every mapping sending $\Gamma$ to
a bounded subset of the odd part of a
locally convex unital graded commutative
algebra $M$ extends to a continuous unital graded algebra
homomorphism $\Lambda\to M$.
The algebra $\Lambda$ is nothing more than $\bigwedge l_1(\Gamma)$,
the free locally convex graded commutative unital algebra
over the Banach space $l_1(\Gamma)$.
\endproclaim

\demo{Proof} Reduces to verification of the universality property
of the kind described in the Theorem for the algebra
$\bigwedge l_1(\Gamma)$, which, in turn,
follows from Theorem 2.1 and the
following statement known in functional analysis:
the topology on the space $l_1(\Gamma)$ is the strongest
complete locally
convex topology making $\Gamma$ into a bounded subset.
\qed\enddemo

We will be interested in the algebra $\G$ over the set
$\Gamma=\N$ of countably many free odd generators.
In a sense, it is one of the most natural
``grassmannian'' algebras.

Denote by $M_n$ the totality of all $n$-element subsets of $\N$
inheriting from $\N$ the natural order, and let
$M=\cup_{n\in\N}M_n$.  For any
$\mu=(\mu_1,\mu_2,\dots,\mu_k)\in M$ put $l(\mu)=k$ and
$\xi^{\mu}:= \xi_{\mu_1}\xi_{\mu_2}\dots\xi_{\mu_k}$.
Then every element $x\in\bigwedge^nl_1$ is uniquely represented as
$$x=\sum_{\mu \in M_n}a_{\mu}\xi^{\mu}$$
with $a_{\mu}\in \C$  and
$$\sum_{\mu \in M_n}\vert a_{\mu}\vert <+\infty$$
An arbitrary element $x\in\G$ can be represented as
a finite sum $x=x_1+\dots+x_n$, where $x_i\in
\bigwedge^{n_i}l_1$
for some $n_i\in\N$.
Every subspace $\bigwedge^nl_1$ is a complete normable
locally convex space, although it does not carry any
distinguished norm.
Actually, the algebra $\G$ can be ``manufactured'' from a
well-known Banach-Grassmann algebra $B_\infty$
\cite{7} by
decomposing the latter algebra into the $l_1$ type sum
of its normed subspaces isomorphic to
$\bigwedge^nl_1$ each, and then rearranging the parts,
this time
as a locally convex direct sum. Recall that $B_\infty$
is the completion of the Grassmann algebra on an infinite
set $\xi_1, \dots, \xi_n, \dots$ of odd generators, endowed
with the norm
$$\norm x \equiv \norm{\sum_{\mu \in M}a_{\mu}\xi^{\mu}}
\overset def\to= \sum_{\mu \in M}\vert a_{\mu}\vert_\C,$$
and an element of $B_\infty$ is any sum
$\sum_{\mu \in M}a_{\mu}\xi^{\mu}$ such that the above norm is
finite.

The algebra $\G$ is, as a locally convex space,
an {\it (LB)}-space (the direct locally convex limit of
the Banach spaces $\oplus_{i=0}^n \bigwedge^nl_1$).

Remark that there is a canonical continuous algebra
monomorphism $j:\G\to B_\infty$. Its restriction to
every subspace $\oplus_{i=0}^n \bigwedge^nl_1$ is an embedding of
locally convex spaces.

\smallpagebreak
\heading
\S 3. Superspaces
\endheading
\smallpagebreak

A {\it locally ringed superspace} \cite{15} is defined by an analogy
with a locally ringed (or geometric) space \cite{4}. It
is a pair $M = (M_0, \Cal S_M)$ consisting of a topological space
$M_0$ and a sheaf $\Cal S_M$ of graded commutative rings on it
such that every stalk $\Cal S_{M,x},~x \in M_0$ is a local ring.
Denote by $m_x$ the radical
of the ring $\Cal S_{M,x}$.
A {\it morphism} between two superspaces $M$ and $N$ is a pair
$\phi = (\phi_0, \phi^\sharp )$ where the mapping
$\phi_0\colon M_0 \to N_0$
is continuous and
$\phi^\sharp  \colon \Cal S_N \to \phi_0^{\star}\Cal S_M$
is a sheaf morphism such that for every $x \in M_0,
\phi^\sharp m_{\phi_0x} \subset m_x$.

This concept is crucial for
supergeometry, because a {\it supermanifold} is
nothing more than a locally ringed superspace
obtained by patching together
the model superspaces of a special form, the {\it superdomains}.
However, what we need here, is the following
particular concept:
a {\it purely odd} supermanifold $\C^{0,q}$
of dimension $(0,q)$ is a superspace
whose underlying topological space is one-pointed and the
algebra of global sections of the (constant) structure sheaf
is the Grassmann algebra $\bigwedge (q)$.

We are interested in the behaviour of the purely odd
superspace $\C^{0,q}$ as $q\to\infty$.
To do that, we need a superspace of infinite purely
odd dimension. There are several different approaches to
infinite-dimensional supermanifold theory \cite{18, 28},
which is far from being in its final form yet.
Happily enough, it would be consistent with all of them, to
assume that a superspace
whose underlying topological space is one-pointed and the
algebra of global sections of the constant structure sheaf
is $\G$, is an example of a supermanifold of dimension
$(0,\infty)$. We will denote this superspace
by $\C^{0,\infty}$.

Our main idea is to consider, within a nonstandard model of
analysis, the superspace $\C^{0,\infty}$,
and to assume that the global function
algebra on an associated superspace $\X$ arising as
the limit of superspaces
$\C^{0,q},~q\to\infty$, or the ``collective effect''
in presence of infinitely many odd coordinates
$\xi_1, \dots ,\xi_n, \dots$ on $\C^{0,\infty}$,
is the {\it nonstandard hull}
of the global function
algebra $\G$ on $\C^{0,\infty}$, that is, the quotient space of
the algebra of all finite elements of $\G$ by the ideal
of all infinitesimals. In a sense, the nonstandard hull is
the ``observable part'' of the algebra $\G$, or (in the spirit of
the French school \cite{11}) its ``shadow'' in the standard world.

\smallpagebreak
\heading
\S 4. Nonstandard hull of the grassmannian algebra $\G$
\endheading
\smallpagebreak

Let $\frak M$ be a higher order set-theoretic structure, and
$^\ast\frak M$ be a higher order nonstandard model of analysis
enlargement of $\frak M$
\cite{26, 30} which is at least $\aleph_1$-saturated.
Let $E\in\frak M$ be a standard locally convex space.
Denote
$$fin~E \overset def\to=
 \{ x\in ^\ast E : p(x)~\text{is a finite element of}~
^\ast \R $$
$$\text{for every
standard continuous prenorm $p$}\},$$
and
$$\mu_E (0) \overset def\to=
 \{ x\in ^\ast E: p(x) \approx 0~\text{for every
standard continuous prenorm $p$}\}.$$
The quotient linear space $\hat E$ of the {\it principal galaxy}
of $E$, $fin~E$, by the {\it monad of zero,} $\mu_E (0)$,
becomes a standard locally convex space if being endowed with a
family of prenorms $\tilde p$, where $p$ runs over the
collection of all continuous prenorms
on $E$,
by letting $\tilde p(y) \overset def\to= st~p(x)$ for any element
$y\in\hat E$ of the form $y=\pi x,~x\in fin~E$, where
$\pi_E:fin~E\to\hat E$ is the quotient linear mapping.
The locally convex space
$\hat E$ is termed the {\it nonstandard hull} of $E$.
The LCS $E$ canonically embeds into $\hat E$ as
a locally convex subspace.
If $\frak M$ is saturated enough, then $\hat E$ is complete.
For more on this, including nonstandard hulls of
internal rather than standard spaces, see \cite{12, 30, 6}.

If $A$ is a standard
locally convex topological algebra then the principal galaxy $fin~A$ is
a (generally speaking, external) subalgebra of $A$, the monad of zero
$\mu_A(0)$ is an ideal of $fin~A$, and the nonstandard hull
$\hat A$ is a locally convex topological algebra belonging to a
variety of algebras generated by $A$.
The algebra $A$ itself is isomorphic to a topological subalgebra of
$\hat A$ in a canonical way.

As a consequence, the nonstandard hull $\H$ of the grassmannian
algebra $\G$ is a complete locally convex
graded commutative algebra, and the nonstandard hull
$\hat B_\infty$ is a graded commutative Banach algebra.
Because of the functorial properties of the operation
of forming nonstandard hull,
the mapping $j: \G\to B_\infty$ gives rise to
a continuous algebra homomorphism $\hat j:\H\to\hat B_\infty$.

\proclaim{4.1. Theorem}
The algebra $\H$ is a locally convex sirect sum
of its complete normable
locally convex subspaces $\widehat{\bigwedge^nl_1}$;
in particular, $\H$ is an (LB)-space.
The mapping $\hat j:\H\to\hat B_\infty$ is a continuous
algebra monomorphism, and
for every $n\in\N$, the restriction of $\hat j$ to
the nonstandard hull $\widehat{\bigwedge^nl_1}$
is an isomorphism of locally convex spaces.
\endproclaim

\demo{Proof}
The algebra $\G$ can be represented as the box product
of spaces $\bigwedge^n l_1$, and it is easy to show that the
nonstandard hull of the  box product of
a standard countable family of LCS's is canonically isomorphic to
the box product of their nonstandard hulls.
The remaining properties are straightforward.
\qed\enddemo

Remark that the nonstandard hull $\hat B_\infty$ is bigger
than just the $l_1$-type sum
of its Banach subspaces $\widehat{\bigwedge^nl_1}$ \cite{22}.

In our paper \cite{22}, we investigated the structure
of the nonstandard hull $\widehat{\bigwedge (\nu)}$
of the (internal)
$\ast$finite dimensional Grassmann algebra
with $\nu$ odd generators, endowed with an $l_1$-type norm.
Although this latter algebra is only a normed subalgebra of
$B_\infty$, and therefore the hull $\widehat{\bigwedge (\nu)}$ is
isometric to a proper closed subalgebra of the hull
$\hat B_\infty$, it is easy to see that
$\widehat{\bigwedge (\nu)}$ is a retract of $\hat B_\infty$,
and therefore
virtually all results can be extended
from the case of $\widehat{\bigwedge (\nu)}$ to
$\hat B_\infty$. In particular, the following holds.

\proclaim{4.2. Theorem {\rm (cf. \cite{22})}}
The spectrum $\Sigma (\hat B_\infty)$ of the nonstandard hull
$\hat B_\infty$ is an inseparable compact space.
The number of
connected components of $\Sigma (\hat B_\infty)$  is uncountable;
however, $\Sigma (\hat B_\infty)$  contains a connected inseparable
subspace. The space  $\Sigma (\hat B_\infty)$
contains a topological copy
of the cube $I^n$ for each natural number $n$,
therefore the topological dimension of  $\Sigma (\hat B_\infty)$
is infinite in any sense.
\qed\endproclaim

The Gelfand spectrum $\Sigma (A)$ of a topological algebra
$A$ is defined as in case of Banach algebras, but now one should
require that $\Sigma (A)$ be formed by {\it continuous}
characters of the algebra $A$ only.
The mapping $\hat j:\H\to\hat B_\infty$ leads to the dual
continuous mapping between Gelfand spectra of both algebras,
$j^\sharp:\Sigma (\hat B_\infty)\to \Sigma (\H)$.
One can show that it is neither one-to-one nor onto.
However, the constructions in the paper \cite{22} were based on the
existence of the following element of $\widehat{\bigwedge (\nu)}$
which is not quasinilpotent and, moreover, has spectral radius $1$:

$$\theta=\pi[{2\over\nu }\sum_{i=1}^{[\frac\nu 2]}\xi_{2i-1}\xi_{2i}]  $$
This element belongs to $\widehat{\bigwedge^nl_1}$ and therefore
is to be found in $\H$ as well. It means that the image
under $j^\sharp$ of any character $\chi$ sending $\theta$ to $1$
is nontrivial. One can verify that the restriction of the
mapping  $j^\sharp$ to any finite-dimensional cube
$I^n$ embedded into $\Sigma (\hat B_\infty)$ as
in our paper \cite{22} is one-to-one, and this way one comes to the
following result about the spectrum of the algebra $\H$.

\proclaim{4.3. Theorem}
The spectrum $\Sigma (\H)$ of the nonstandard hull
$\H$ is an inseparable compact space.
The number of
connected components of $\Sigma (\H)$  is uncountable;
however, $\Sigma (\H)$  contains a connected inseparable
subspace. The space  $\Sigma (\H)$
contains a topological copy
of the cube $I^n$ for each natural number $n$,
therefore the topological dimension of  $\Sigma (\H)$
is infinite in any sense.
\qed\endproclaim

\smallpagebreak
\heading
\S 5. Automorphisms of the algebra $\H$
\endheading
\smallpagebreak

Let $u\in GL (l_1)$ be an automorphism of the Banach space
$l_1$. In view of existence of canonical embedding
$l_1 \hookrightarrow \bigwedge l_1$,
the operator $u$ becomes a linear mapping $u: l_1 \to \G$.
Since the set $u(\Xi)$ is bounded, the operator $u$ can be extended
to a standard automorphism $\tilde u$ of $\G$,
which in turn
is lifted to an automorphism, $\hat u$, of
the hull algebra $\H$. This way the general linear group
$GL (l_1)$ acts on the algebra $\H$ by automorphisms.
We denote this action by $\tau$ and endow the group
$GL (l_1)$ with the uniform operator topology.

\proclaim{5.1. Theorem}
The action $\tau$ is analytic as a mapping
$GL (l_1)\times \H \to \H$.
\endproclaim

\demo{Proof}
Since every subspace $\widehat{\bigwedge^nl_1}$ of $\H$
is closed under the action $\tau$, the action itself
decomposes into the locally convex
direct sum of actions $\tau_n$ of
$GL (l_1)$ on $\widehat{\bigwedge^nl_1}$,
and it is sufficient to prove the analyticity of each mapping
$\tau_n:GL (l_1)\times \widehat{\bigwedge^nl_1}
\to \widehat{\bigwedge^nl_1}$, as $n\in\N$.
Since $GL (l_1)$ is a Banach-Lie group and the locally convex space
$\widehat{\bigwedge^nl_1}$ is complete normable,
the analyticity of the action
would follow at once from its continuity as a homomorphism
$GL (l_1)\to GL(\widehat{\bigwedge^nl_1})$ with respect to the
uniform operator topology on the latter group.
It is enough to check the continuity of this homomorphism
at the identity element only. Let us endow, for the purpose of
our proof, the locally convex space
$\bigwedge^n~^\ast l_1$ with a norm induced from the Banach algebra
$^\ast B_\infty$
(a particular choice of the norm does not affect
the uniform operator topology).
Let an $\e>0$ be arbitrary standard with $\e<1$;
we will find a neighbourhood of
identity, $U$, in $GL (l_1)$ such that for every $g\in U$ and
$x\in \bigwedge^n~^\ast l_1$ with $\norm x \leq 1$ one has
$\norm{\tau_g(x)-x}\leq\e$.
Choose, using the uniform continuity of the action of
$GL (l_1)$ on $l_1$, a neighbourhood of identity,
$U\subset GL (l_1)$, such that
for every $g\in ^\ast U$ and
$x\in ^\ast l_1$ with $\norm x \leq 1$ one has
$\norm{\tau_g(x)-x}\leq\e/2^{n-1}n$.
Since the unit ball in
$\bigwedge^n~^\ast l_1$ is the convex circled envelope
of the set of monomials $\xi^\mu$, $\mu\in ^\ast M_n$,
one can assume that in the statement we wish to prove,
$x=\xi^\mu,~\mu\in ^\ast M_n$,
and, even more, that
$x=\xi^\mu,~\mu\in M_n$ (by $\ast$-Transfer). Let $g\in U$ and
$x=\xi_{\mu_1}\dots \xi_{\mu_n}$.
One has
$$\norm{\tau_g(x)-x}\equiv
\norm{\tau_g(\xi_{\mu_1})\dots \tau_g(\xi_{\mu_n})-
\xi_{\mu_1}\dots \xi_{\mu_n}} $$
$$\leq \norm{\tau_g(\xi_{\mu_1})\tau_g(\xi_{\mu_2})\tau_g(\xi_{\mu_3})
\dots \tau_g(\xi_{\mu_n})-
\xi_{\mu_1}\tau_g(\xi_{\mu_2})\tau_g(\xi_{\mu_3})
\dots \tau_g(\xi_{\mu_n})} $$
$$+\norm{\xi_{\mu_1}\tau_g(\xi_{\mu_2})\tau_g(\xi_{\mu_3})
\dots \tau_g(\xi_{\mu_n})-
\xi_{\mu_1}\xi_{\mu_2}\tau_g(\xi_{\mu_3})
\dots \tau_g(\xi_{\mu_n})}+\dots $$
$$+\norm{\xi_{\mu_1}\xi_{\mu_2}\xi_{\mu_3}\dots
\xi_{\mu_{n-1}}\tau_g(\xi_{\mu_n})-
\xi_{\mu_1}\xi_{\mu_2}\xi_{\mu_3}\dots \xi_{\mu_{n-1}}\xi_{\mu_n}} $$
$$ \leq \norm{[\tau_g(\xi_{\mu_1})-\xi_{\mu_1}]
\tau_g(\xi_{\mu_2})\tau_g(\xi_{\mu_3})
\dots \tau_g(\xi_{\mu_n})}$$
$$+\norm{\xi_{\mu_1}[\tau_g(\xi_{\mu_2})-\xi_{\mu_2}] \tau_g(\xi_{\mu_3})
\dots \tau_g(\xi_{\mu_n})} $$
$$+\dots + \norm{\xi_{\mu_1}\xi_{\mu_2}\xi_{\mu_3}\dots \xi_{\mu_{n-1}}
[\tau_g(\xi_{\mu_n})-\xi_{\mu_n}]}$$
$$\leq n\cdot 2^{n-1}\cdot\frac\e {2^{n-1}n} = \e,$$
as desired.
\qed\enddemo

\smallpagebreak
\heading
\S 6. Analyticity in the spectrum of $\H$
\endheading
\smallpagebreak

\proclaim{6.1. Theorem}
There exists a homeomorphism $\Phi$ from the unit disc $D\subset \C$
into $\s$ such that for every element
$a\in\H$ the composition $\hat a\circ\Phi$ is a holomorphic
function from $D$ to $\C$.
\endproclaim

\demo{Proof}
For every $z\in D$, the dilatation of  $l_1$ by the element
$e^z$, i.e., the mapping $x\mapsto e^zx$, belongs to
$GL(l_1)$; we will denote this dilatation simply by $e^z$.
Now fix a continuous character $\chi$
(that is, an element of $\s$) sending the element
$\theta$ (constructed in $\S 4$) to $1$.
(Such a character exists in the spectrum of the Banach algebra
$\hat B_\infty$, as it is well known, and it is enough to take
its restriction to $\H$ as $\chi$.)
For each $z\in D$ define a new character $\chi_z$ by letting
$\chi_z(a)=\chi (\tau_{e^z}a)$ for all $a\in\H$.
The arising mapping $\Phi : z\mapsto \chi_z$ from $D$ to
$\s$ is one-to-one (because $\chi_z(\theta)=e^{2z}$ for all
$z\in D$) and is continuous. Since the topology on $\s$ is that
of pointwise convergence on elements of $\H$, then
in order to prove the continuity of $\Phi$,
it is enough to show
that all mappings $z\mapsto \chi_z (a)$ from $D$ to $\C$
are continuous as $a\in\H$. This follows from a stronger result:
for every $a\in\H$, the mapping $z\mapsto \chi_z (a)$
(in other form, $\hat a\circ\Phi$)
is analytic, as the composition of a chain of
analytic maps between complete locally convex spaces:

$$
\CD
@. D @>z\mapsto e^z>> GL(l_1) @. @.\\
D @>z\mapsto (z,a)>> \times @. \times @>\tau>> \H @>\chi>> \C \\
@. \{a\} @>>>\H @. @.
\endCD
$$
\qed\enddemo

\smallpagebreak
\heading
\S 7. Conclusion: further developments
\endheading
\smallpagebreak

Since we are working in an analytic rather than purely algebraic
context,
it seems reasonable to modify the definition of
a locally ringed ($=$ geometric) superspace by letting
algebras of
superfunctions ($=$ sections of the strucutre sheaf) carry a
topology, which, for a number of
reasons, should be complete and
locally {\it m}-convex in the sense of
\cite{1, 16}. The restriction homomorphisms
must be, of course, continuous.
One can show that, like in the case of algebraic prime spectra
\cite{4}, there exists a universal solution to the problem
of associating to any graded commutative Banach algebra $A$
a superspace in the above sense;
the underlying topological space
of the universal superspace is the Gelfand spectrum
$\Sigma (A)$. (A related construction --- in the ``purely
even'' case --- can be found in \cite{25}.)
The  superspace $\X$ which we consider
as the ``collective effect'' arising in
$\C^{0,\infty}$
is just the solution to a universal problem of such
kind for the hull algebra $\H$.
(A superspace
denoted by the same letter $\X$ in \cite{19, 22}
is some proper subsuperspace.)
The reduced subsuperspace, $\X_{red}$,
of $\X$ \cite{15} is obtained by taking quotient of the structure
sheaf of $\X$ by the ideal of nilpotents (or, in our topologized case,
quasinilpotents).

It turns out that
the non-nilpotent derivations of the hull
algebra $\H$ are abundant, and
they give birth to vector fields
on the geometric space $\X_{red}$.
(One can obtain those derivations as elements of
the nonstandard hull
of the well-known \cite{8} Lie superalgebra of derivations of
a Grassmann algebra.)
Thereby, the space $\X_{red}$
admits a substantial
differential geometry of its own.

Another possibility to treat
the problem is to invoke
the synthetic approach, based on topos theory,
and the first step may be to give an appropriate locally convex
nonstandard hull of an infinite-dimensional Grassmann algebra
the structure of a smooth algebra in the sense of \cite{17}.

\smallpagebreak
\heading Acknowledgments
\endheading
\smallpagebreak

This investigation
\footnote{Partially supported
by a VUW 1992
Small Research Grant V212/451/RGNT/594/153.}
was conceived on the shores of the
Adriatic Sea in December 1990,
and performed on the shores of the Tasman Sea almost two years later.
The author is thankful to
Gianni Landi and Cesaro Reina for inviting him to come to
SISSA (Trieste, Italy) for
a week, and to
Rob Goldblatt for inviting him to come to VUW
(Wellington, NZ) for good.

\bye